%% file: tot.tex
\documentclass{article}
\usepackage{amssymb}
\usepackage{latexsym}
\newcommand{\N}{{\mathbb N}}

\newcommand{\Z}{{\mathbb Z}}

\newcommand{\bx}{\hfill{$\Box $ }}

\newcommand{\Img}{\textrm{Im}\;}

\newcommand{\mod}{\textrm{mod}\;}

\newtheorem{corollary}{Corollary}[section]
\newtheorem{proposition}{Proposition}[section]
\newtheorem{theorem}{Theorem}[section]
\title{On the image of Euler's totient function}
\author{R.Coleman\\Laboratoire LJK-Universit\'e de Grenoble,\\ Tour-IRMA,\\51, rue des Math\'ematiques,\\Domaine universitaire de Saint-Martin-d'H\`eres, France.}
\begin{document}
\maketitle

\noindent Euler's totient function $\phi$ is the function defined on the positive natural numbers
$\N^*$ in the following way: if $n\in \N^*$, then $\phi (n)$ is the cardinal of
the set 
$$
\{x\in \N^*:1\leq x\leq n, (x,n)=1\},
$$ where $(x,n)$ is the pgcd of $x$ and $n$. Thus $\phi (1)=1$, $\phi (2)=1$, 
$\phi (3)=2$, $\phi (4)=2$, and so on. The principle aim of this article is to 
study certain aspects of the image of the function $\phi$.\\

\input{elprop}

\input{bounds}

\input{2p}

\input{set}
\input{struct}

\input{concl}

To close I would like to thank Dominique Duval and Mohamed El Methni for reading the text
 and offering helpful suggestions.\\
 
\newpage

\input{bibliogr}
%\printindex

\end{document}

%% file: elprop.tex
\section{Elementary properties} Clearly $\phi (p)=p-1$, for any prime number 
$p$ and, more generally, if $\alpha \in \N^*$, then $\phi (p^{\alpha})=p^{\alpha }-
p^{\alpha -1}$. This follows from the fact that the only numbers which are not 
coprime with $p^{\alpha}$ are multiples of $p$ and
there are $p^{\alpha -1}$ such multiples $x$ with $1\leq x\leq p^{\alpha}$.\\

It is well-known that $\phi$ is multiplicative, i.e., if $m$ and $n$ are 
coprime, then $\phi (mn) = \phi (m)\phi (n)$. If $n\geq 3$ and the prime 
decomposition of $n$ is 
$$
n = p_1^{\alpha _1}\ldots  p_s^{\alpha _s},
$$
then from what we have seen 
$$
\phi (n) = \prod _{i=1}^s(p_i^{\alpha _i}-p^{\alpha _i-1}) = 
n\prod _{i=1}^s(1-\frac{1}{p_i}).
$$
Notice that 
$$
p^{\alpha }- p^{\alpha -1}=p^{\alpha -1}(p-1).
$$
This implies that, if $p$ is odd or $p=2$ and $\alpha >1$, then 
$p^{\alpha }- p^{\alpha -1}$ is even. Hence, for $n\geq 3$, $\phi (n)$ is even.
Thus the image of $\phi$ is composed of the number $1$ and even numbers.\\

The following property is simple but very useful.

\begin{proposition}\label{prop1} If $p,m\in \N^*$, with $p$ prime, and $n=pm$, then $\phi
(n)=(p-1)\phi (m)$, if $(p,m)=1$, and $\phi (n)=p\phi (m)$, if $(p,m)\neq 1$.
\end{proposition}

\noindent \textsc{proof} If $(p,m)=1$, then we have
$$
\phi (n) =\phi (p)\phi (m) = (p-1)\phi (m).
$$
Now suppose that $(p,m)\neq 1$. We may write $m=p^{\alpha}m'$, with 
$\alpha \geq 1$ and $(p,m')=1$. Thus 
$$
\phi (n) = \phi (p^{\alpha +1})\phi (m') = p^{\alpha}(p-1)\phi (m').
$$
However 
$$
\phi (m) = p^{\alpha -1}(p-1)\phi (m')
$$
and so $\phi (n)= p\phi (m)$.\bx

\begin{corollary} \label{cor1a} $\phi (2m)=\phi (m)$, if and
only if $m$ is odd.
\end{corollary}

\noindent \textsc{proof} If $m$ is odd, then $(2,m)=1$ and so $\phi
(2m)=(2-1)\phi (m)=\phi (m)$. If $m$ is even, then $(2,m)\neq 1$, hence $\phi
(2m) = 2\phi (m)\neq \phi (m)$. This ends the proof.\bx

%% file: bounds.tex
\section{Bounds on $\phi ^{-1}(m)$} Let $m\in \N^*$ and consider the inverse
image $\phi ^{-1}(m)$ of $m$, i.e.
$$
\phi ^{-1}(m) = \{n\in \N^*:\phi (n)=m\}.
$$
There are several questions we might ask. First, is the set $\phi ^{-1}(m)$
empty and, if not is it finite. It is easy to see that $\phi ^{-1}(1)=\{1,2\}$ 
and as we have already seen, $\phi ^{-1}(m)$ is empty if $m$ is odd. It remains 
to consider the case where $m$ is an even number. The following result, due to 
H. Gupta \cite{gupta}, helps us to answer these questions.
 
\begin{proposition} Suppose that $m$ is an even number and let us set
$$
A(m) = m\prod _{p-1|m}\frac{p}{p-1},
$$
where $p$ is prime. If $n\in \phi ^{-1}(m)$, then $m<n\leq A(m)$. 
\end{proposition}

\noindent \textsc{proof} Clearly, if $\phi (n)=m$, then $m<n$. On the other
hand, if $\phi (n) =m$ and $n = p_1^{\alpha _1}\ldots  p_s^{\alpha _s}$, then
$$
m = n\prod _{i=1}^s\frac{p_i-1}{p_i} \Longrightarrow 
	n=m\prod _{i=1}^s\frac{p_i}{p_i-1}.
$$
However, if $p|n$, then from the first section we know that $p-1|\phi (n)$ and it
follows that, for each $p_i$, $p_i-1|m$. Hence $n\leq A(m)$.\bx\\

The proposition shows that the inverse image of an element $m\in \N^*$ is always
finite. It also enables us to determine whether a given number
$m$ is in the image of $\phi$: we only need to determine $A(m)$, and then
calculate $\phi (n)$ for all integers $n$ in the interval $(m,A(m)]$.\\

\noindent {\bf Examples} {\bf 1.} The divisors of $4$ are 1,2 and 4. Adding 1 
to each of these numbers we obtain 2, 3 and 5, all of which are prime numbers. 
Thus 
$A(4) = 4\cdot\frac{2}{1}\cdot\frac{3}{2}\cdot\frac{5}{4} = 15$.
To find the inverse image of 4, it is sufficient to consider numbers between 5
and 15. In fact, $\phi ^{-1}(4)=\{5, 8, 10, 12\}$.\\
{\bf 2.} The divisors of 14 are 1, 2, 7 and 14. However, if we add 1 to each of 
these numbers we only find a prime number in the first two cases. Thus
$A(14)=14\cdot\frac{2}{1}\cdot\frac{3}{2}=42$. If we consider the numbers $n$ between 15
and 42, we find $\phi (n)\neq 14$ and so $\phi ^{-1}(14)=\emptyset$.\\

\noindent {\bf Remark} The example $m=14$ shows that there are even
numbers which are not in the image of $\phi$.\\

$A$ is a function defined on $\{1\}\cup 2N^*$. Let us look at the first values
of $A$ with the corresponding values of $\phi$ (when defined):
\begin{center}
\begin{tabular}{|c|c|c|}\hline
$m$ & $A(m)$ & $\phi (A(m))$\\ \hline
$1$ & $2$    & $1$ \\ \hline	 
$2$ & $6$    & $2$ \\ \hline
$4$ & $15$    & $8$ \\ \hline
$6$ & $21$    & $12$ \\ \hline
$8$ & $30$    & $8$ \\ \hline
$10$ & $33$    & $20$ \\ \hline
$12$ & $\frac{455}{8}$    & $-$ \\ \hline
$14$ & $42$    & $20$ \\ \hline
\end{tabular}
\end{center}
\vspace{3mm}
From the table we see that $A(m)$ may be odd, even or a fraction and that $A$ is
not increasing ($A(12)>A(14)$). Also, we may have $\phi (A(m))=m$ or 
$\phi (A(m))>m$.\\

It is interesting to consider the special case $m=2^k$. The only divisors of
$2^k$ are $1,2,2^2 \ldots ,2^k$. If we add 1 to each of these numbers, we obtain
numbers of the form $2^i+1$, where $0\leq i\leq k$. For $i\geq 1$, such
numbers are prime only if $i$ is a power of $2$. A number of the form
$F_n=2^{2^n}+1$ is said to be a Fermat number. Therefore
$$
A(2^k) = 2^k\cdot 2\cdot \prod \frac{F_n}{F_n-1}=
	2^{k+1}\cdot\prod \frac{F_n}{F_n-1},
$$
where the product is taken over Fermat numbers $F_n$ which are prime and such 
that $F_n-1|2^k$. For example, if $m=2^5=32$, then $F_0$, $F_1$ and $F_2$ are
the only Fermat numbers $F_n$ such that $F_n-1|m$. In addition these Fermat
numbers are all prime. Therefore 
$
A(m) = 64\cdot \frac{3}{2}\frac{5}{4}\cdot \frac{17}{16} = \frac{255}{2}
$. The Fermat numbers $F_0, \ldots , F_4$ are all prime; however $F_5,\ldots ,F_{12}$ are 
composite numbers and it has been shown that for many
other numbers $n$, $F_n$ is composite. In fact, up till now no Fermat number
$F_n$ with $n\geq 5$ has been found to be prime. If there are no prime Fermat
numbers with $n\geq 5$, then for $2^k\geq 2^{16}$ we have 
$$
A(2^k)=2^{k+1}\cdot \frac{F_0}{2}\cdot\frac{F_1}{2^2}\cdot\frac{F_2}{2^4}
	\cdot\frac{F_3}{2^8}\cdot\frac{F_4}{2^{16}}=2^{k-30}F_0F_1F_2F_3F_4, 
$$
which is an integer for $k\geq 30$.\\

Before closing this section, let us consider the upper bound on odd elements of
$\phi ^{-1}(m)$. From Corollary \ref{cor1a}, we know that if $n$ is odd and $n\in
\phi ^{-1}(m)$, then $2n\in \phi ^{-1}(m)$, therefore an upper bound on odd
elements of $\phi ^{-1}(m)$ is $\frac{A(m)}{2}$. We should also notice that at
least half of the elements in $\phi ^{-1}(m)$ are even. In fact, $\phi ^{-1}(m)$
may be non-empty and contain very few odd numbers, or even none. For example, 
the only odd number in $\phi ^{-1}(8)$ is $15$ and $\phi ^{-1}(2^{32})$ contains
only even numbers (see Theorem \ref{th5a} further on).

%% file: 2p.tex
\section{The case $m=2p$} We now consider in some detail the case where $p$ is
prime and $m=2p$.

\begin{theorem} \label{thm3a} If $p$ is a prime number, then $2p$ lies in the image of
$\phi$ if and only if $2p+1$ is prime.
\end{theorem} 

\noindent \textsc{proof} If $2p+1$ is prime, then $\phi (2p+1)=2p$ and so 
$2p\in \Img\phi$.\\

\noindent Suppose now that $2p\in \Img\phi$. If $p=2$, then $2p=4\in \Img\phi$,
car $\phi (5)=4$ and $2p+1=5$, which is prime. Now suppose that $p$ is an odd
prime. As $2p\in\Img \phi$, there exists $n$ such that $\phi (n)=2p$. If $n=2^k$, 
then $2^{k-1}=\phi (2^k)=2p$, which is clearly impossible, because $p$ is an 
odd number greater than 1. Hence there is an odd prime $q$ such that $n=qs$. 
There are two cases to consider: 1. $q\not |s$, 2. $q|s$. We will handle each 
of these cases in turn, using Proposition \ref{prop1}.\\   

\noindent Case 1. We have $2p=(q-1)\phi (s)$ which implies that $q-1|2p$. The 
only divisors of $2p$ are $1$, $2$, $p$ and $2p$. As $q\neq 2$, the possible 
values for $q-1$ are 2, $p$ or $2p$ and hence for $q$ are 3, $p+1$ or $2p+1$. However,
$p+1$ is not possible, because $p+1$ is even and hence not prime. If $q=3$, then
$\phi (s)=p$. As $p$ is an odd number greater than 1, this is not possible. It
follows that $q=2p+1$ and so $2p+1$ is prime.\\

\noindent Case 2. Here we have $2p=q\phi (s)$ and so $q|2p$. The possible 
values of $q$ are 1, 2 , $p$ or $2p$. However, as $q$ is an odd prime, we must 
have $q=p$. This implies that $\phi (s)=2$. Also, $q|s$ and so $q-1|\phi (s)$. 
It follows that $p=q=3$ and hence $2p+1=7$, a prime number.\\

\noindent As in both cases $2p+1$ is prime, we have proved the result.\bx\\

\noindent {\bf Definition} A prime number of the form $2p+1$, with $p$ prime, is
said to be a safe prime. In this case, the prime number $p$ is said to be a 
Sophie Germain prime.\\

\noindent {\bf Remarks 1.} The theorem does not generalize to odd numbers. Certainly,
if $s$ is odd and $2s+1$ is prime, then $2s\in \Img \phi$. However, it may be so
that $2s+1$ is not prime and $2s\in \Img\phi$. For example, $54=\phi (81)$ and
so $54\in\Img\phi$. However, $54=2\cdot 27$ and $2\cdot 27+1=55$, which is not
prime.\\
\noindent {\bf 2.} Let us now consider $\phi ^{-1}(2p)$. From Theorem 
\ref{thm3a}, the set $\phi ^{-1}(2p)$ is empty if $2p+1$ is not prime. If 
$2p+1$ is prime, then we have
$$
A(2p) = 2p\cdot \frac{2}{1}\cdot\frac{3}{2}\cdot\frac{2p+1}{2p} = 6p+3.
$$
Hence, if $n\in \phi ^{-1}(2p)$, then $2p<n\leq 6p+3$. It is worth noticing that
$$
\phi (6p+3) = \phi (3)\phi (2p+1) = 2\cdot 2p = 4p.
$$

\begin{corollary} If $2p\in \Img \phi$, then $2^kp\in \Img \phi$, for $k\geq 1$.
\end{corollary} 

\noindent \textsc{proof} For $k=1$ the result is already proved, so suppose that
$k\geq 2$. As $2p+1$ is an odd prime, $2^k$ and $2p+1$ are coprime. Therefore
$$
\phi (2^k(2p+1))=\phi (2^k)\phi (2p+1)=2^{k-1}2p=2^kp.
$$
This ends the proof.\bx\\

\noindent {\bf Remarks} {\bf 1.} From the corollary we deduce that, for every prime
$p$ such that $2p+1$ is prime, there is an infinity of numbers $a$ such that 
$$
a\equiv 0 (\mod p)
$$
and $a\in \Img \phi$.\\
\noindent {\bf 2.} If $2p\notin \Img \phi$, then we cannot say that 
$2^kp \notin \Img \phi$ for any $k\geq 2$. For example, $14\notin \Img \phi$,
but $28\in \Img \phi $.\\

%% file: set.tex
\section{Sets of elements of the form $2p$} In this section we will need an
important theorem due to Dirichlet: 

\begin{theorem} If $n\in \N^*$ and $(a,n)=1$, then there is an infinite number
of prime numbers $p$ such that 
$$
p\equiv a (\mod n ).
$$
\end{theorem}
 
\noindent Chapman \cite{chapman} has recently given a relatively elementary proof of 
this result.\\

Let us consider the prime numbers $p$ in the interval $[1,50]$. There are $15$ such
numbers, namely 
$$
2,3,5,7,11,13,17,19,23,29,31,37,41,43,47.
$$ 
For seven of these numbers, $2p+1$ is prime, i.e. $2,3,5,11,23, 29,41$, and for
the others $2p+1$ is not prime. Therefore
$$
4,6,10,22,46,58,82\in \Img \phi \qquad\mathrm{and}\qquad 14,26,34,38,62,74,86,94\notin \Img
\phi .
$$

It is natural to ask whether there is an infinite number of distinct primes $p$
such that $2p\in \Img \phi$ (resp. $2p\notin \Img \phi$). Our question amounts 
to asking whether there is an
infinite number of safe primes, or equivalently an infinite number of Sophie
Germain primes. Up till now this question has not been answered.
The largest known safe prime, found by David Undebakke in January, 2007, is
$48047305725.2^{172404}-1$. We can say a lot more concerning numbers of the 
form $2p$ which are not in the image of $\phi$. 

\begin{theorem} For any odd prime number $p$, there is an infinite set $S(p)$ of 
prime numbers $q$ such that $2q\notin \Img \phi$ and $p|2q+1$.
\end{theorem}

\noindent \textsc{proof} If $p$ is an odd prime, then $\frac{p-1}{2}$ is a
positive integer and $(\frac{p-1}{2},p)=1$. From Dirichlet's theorem, we know
that there is an infinite number of prime numbers of the form
$q=\frac{p-1}{2}+kp$, with $k\in \Z$. Then 
$$
2q+1 = 2\left(\frac{p-1}{2}+kp\right)+1 = p+2kp = p(1+2k).
$$
From Theorem \ref{thm3a}, we know that $2q\notin \Img \phi$ and clearly
$p|2q+1$.\bx\\

\noindent {\bf Remark} The sets $S(p)$ may have common elements, but in general
are distinct. For example, $14$ is in $S(3)$ and
$S(5)$, but not in $S(7)$, $26$ is in $S(3)$, but not in $S(5)$ and $S(7)$ and
$34$ is in $S(5)$ and $S(7)$, but not in $S(3)$.\\

The following result follows directly from the theorem.

\begin{corollary} \label{cor4a} For every odd prime number $p$, there exists an
infinite subset $\tilde{S}(p)$ of $\N^*$ such that
$$
a\equiv -1 (\mod p).
$$
and $a\notin \Img \phi$, when $a\in \tilde{S}(p)$.
\end{corollary}

\noindent {\bf Remark} At least for the moment, we cannot say that there is an infinity of odd prime
numbers $p$ such $2p\in \Img \phi$. However, we can say that that there is an
infinity of odd numbers $s$ such that $2s\in \Img \phi$. Here is a proof. From 
Dirichlet's theorem, there is an infinity of prime numbers $p$ such that 
$p\equiv 3 (\mod 4)$. Thus, for each of these primes, there exists 
$k\in \Z$ such that 
$$
p = 3 + 4k = 1 + 2(1+2k).
$$
As $1 + 2(1+2k)$ is prime $2(1+2k)\in \Img \phi$ and of course $1+2k$ is odd.\\

%% file: struct.tex
\section{Structure of $\phi ^{-1}(m)$} If $n$ is an odd solution of the 
equation $\phi (n)=m$, then $2n$ is also a solution (Corollary \ref{cor1a}). It
follows that the equation can have at most half its solutions odd. It is 
natural to look for cases where there are exactly the same number of odd and 
even numbers of solutions.\\

First let us consider the case where $m=2p$ and $p$ is an odd prime.
If $p=3$, then $6<n\leq 21$. A simple check shows that $\phi
^{-1}(6)=\{7,9,14,18\}$. There are two odd and two even solutions of the
equation $\phi (n)=6$.

\begin{proposition} If $p$ is a prime number such that $p\geq 5$ and $2p+1$ is 
prime, then $\phi ^{-1}(2p)$ contains exactly one odd and one even element,
namely $2p+1$ and $4p+2$.
\end{proposition}

\noindent \textsc{proof} If $\phi (n)=2p$, then, from Remark 2. following 
Theorem \ref{thm3a}, $2p<n\leq 6p+3$. The divisors of
$2p$ are $1$, $2$, $p$ and $2p$ and so the only possible prime divisors of $n$
are $2$, $3$ and $2p+1$. If $n=2p+1$ or $4p+4$, then $\phi (n)=2p$ and there can
be no other multiple $n$ of $2p+1$ such that $\phi (n)=2p$. If $n$ is not a
multiple of $2p+1$ and $\phi (n)=2p$, then $n$ must be of the form
$n=2^{\alpha}\cdot3^{\beta}$. If $\alpha \geq 3$, then $4|2p$, which is
impossible. Also, if $\beta \geq 2$, then $3|2p$. However, this is not possible,
because $p\geq 5$. Therefore $\alpha \leq 2$ and $\beta \leq 1$. As $n>2p\geq
10$, the only possibility is $n=12$. As $\phi (12)=4$, this is also impossible.
The result now follows.\bx\\

We can generalize this result to odd numbers in general. We will write $O(m)$
(resp. $E(m)$) for the set of odd (resp. even) solutions of the equation 
$\phi (n)=m$. 

\begin{theorem} \label{prop4a} If $s$ is odd and $s\geq 3$, then $O(2s)=E(2s)$. 
\end{theorem} 

\noindent \textsc{proof} If the equation $\phi (n)=2s$ has no solution, then 
there is nothing to prove, so suppose that this is not the case. If $n$ is an 
odd solution of the equation, then $2n$ is
also a solution. Hence $O(2s)\leq E(2s)$.\\
If $n$ is an even solution, then 
we may write $n=2^{\alpha}t$, with $\alpha \geq 1$ and $t$ odd. If $t=1$, then 
$\phi (n)$ is a power of $2$, which is not possible. It follows that $t\geq 3$. 
If we now suppose that $\alpha >1$, then $\phi (n)= 2^{\alpha -1}\phi (t)$ and, 
as $\phi (t)$ is even, $4$ divides $\phi (n)=2s$, which is not possible. 
Therefore $\alpha =1$ and so $n=2t$. As $\phi (t)=\phi (2t)=\phi (n)$, we must
have $O(2s)\geq E(2s)$ and the result now follows.\bx

\begin{corollary} There is an infinity of numbers $m$ such that $\phi ^{-1}(m)$
is non-empty and composed of an equal number of odd and even numbers.
\end{corollary}

\noindent \textsc{proof} It is sufficient to recall that there is an infinity 
of primes $p$ such that $p=1+2s$, with $s$ odd (remark after Corollary 
\ref{cor4a}).\bx \\

%\begin{proposition} We have the relation
%$$
%S_e(2m) = 2S_o(2m) \cup 2S_e(m).
%$$
%\end{proposition}

%\noindent \textsc{proof} If $n\in S_o(2m)$, then $\phi (n)=2m$ As $n$ is odd,
%$\phi (2n)=2m$ and so $2n\in S_e(2m)$. Now let us consider the case where 
%$n\in S_e(m)$. We may write $n=2^{\alpha}s$, where $\alpha \geq 1$ and
%$(2,s)=1$. Then 
%$$
%\phi (n) = 2^{\alpha -1}\phi (s) = m
%$$
%and 
%$$
%\phi (2n) = \phi (2^{\alpha +1}s) = 2^{\alpha}\phi (s) =2(2^{\alpha -1}\phi (s))
%=2m
%$$
%and so $2n\in S_e(2m)$. We have shown that $2S_o(2m) \cup 2S_e(m)\subset
%S_e(2m)$.\\
%Now let us look at the other inclusion. Suppose that $n\in S_e(2m)$ and that 
%$n\notin 2S_o(2m)$. Then $n=2^{\alpha}s$, with $s$ odd and $\alpha \geq 2$. We
%have
%$$
%\phi (n) = 2^{\alpha -1}\phi (s) = 2m \Longrightarrow 2^{\alpha -2}\phi (s)=m.
%$$
%Also,
%$$
%\phi (\frac{n}{2}) = \phi (2^{\alpha -1}s) = A^{\alpha -2}\phi (s) =m
%\Longrightarrow \frac{n}{2}\in S_e(m).
%$$
%It follows that $n\in 2S_e(m)$. Therefore 
%$$
%S_e(2m) \subset 2S_o(2m) \cup 2S_e(m)
%$$
%and the result now follows.\bx\\

%\begin{corollary} If $s$ is an odd number with $s\geq 3$ and $m=2s$, then the
%number of solutions of the equation $\phi (n)=m$, has the same number of odd and
%even solutions. 
%\end{corollary}

%\noindent \textsc{proof} It is sufficient to notice that in this case $S_e(m)$
%is empty.\bx\\

%This corollary enables us to extend Proposition \ref{prop4a}. We have already 
%seen that $54\in \Img\phi$, however it is not covered by the previous 
%proposition.\\  

At the other extreme is the case where $\phi ^{-1}(m)$ is non-empty and
composed entirely of even numbers. In considering this question, the following
result due to Gupta \cite{gupta}, is useful. We will give a modified proof of
it.

\begin{theorem} \label{th5a} $O(2^k)=0$ or $O(2^k)=1$. 
\end{theorem}

\noindent \textsc{proof} If $k=0$, then we have $\phi ^{-1}(2)=\{3,4,6\}$ and the result
follows. Suppose now that $k>0$ and that $\phi (n)=2^k$, with $n$ odd. If $p$ 
is an odd prime such that $p-1$ divides $2^k$, then $p$ must be of the form 
$2^i+1$. The only primes of this form are Fermat numbers, hence if
$n=p_0^{\alpha _1}\ldots p_r^{\alpha _r}$ is the decomposition of $n$ as a
product of primes, then each $p_i$ must be a Fermat number. If we allow 
$\alpha _i=0$, then we may suppose that $p_i=F_i$. Thus
$$
2^k = \phi (F_0^{\alpha _0})\ldots \phi (F_r^{\alpha _r})=\prod _{\alpha _i\geq
1}F_i^{\alpha _i-1}(F_i-1).
$$
If $\alpha _i>1$ for some $i$, then the product is not a power of $2$, hence we
must have $\alpha _i=0$ or $\alpha _i=1$, for all $i$, and so 
$$
2^k = \alpha _02^{2^0}\ldots \alpha _r2^{2^r}.
$$
Clearly $\alpha _0\ldots \alpha _r$ is $k$ written in binary form. Therefore 
there can be at most one odd number $n$ such that $\phi (n)=2^k$. If 
$\alpha _i=1$ only when $F_i$ is prime, then
there exists $n$ odd such that $\phi (n)=2^k$ and $s(k)=1$. On the other hand,
if there is an $\alpha _i$ such that $F_i$ is not prime, then there does not
exist an odd number $n$ such that $\phi (n)=2^k$.\bx\\

If $k<32$, then $k$ can be written in binary form as
$\alpha _0\ldots \alpha _4$. As $F_0,\ldots F_4$ are all primes, there exists an
odd number $n$ such that $\phi (n)=2^k$. However, for $2^{32}$, this is not the
case, because $F_{5}$ is not prime. Up till now no Fermat number $F_n$, with
$n\geq 5$, has found to be prime, so it would seem that, for $k\geq 5$, there is
no odd number in the set $\phi ^{-1}(2^k)$. This suggests that there is 
an infinity of numbers $m$ for which $\phi ^{-1}(m)$ is non-empty and only
composed of even numbers.

%% file: concl.tex
\section{Conclusion} We have seen that the inverse image $\phi ^{-1}(m)$ is an empty set for any odd
positive integer $m>1$ and also for an infinite number of even positive
integers. We have also seen that, when $\phi ^{-1}(m)$ is non-empty, its 
cardinal can be odd or even and that it can have at most half of its members
odd; it is also possible that all its members are even. Up till now no number
$m$ has been found such that $\phi ^{-1}(m)$ contains only one element.
Carmichael conjectured that such a case does not exist, but this is yet to be
proved (or disproved). However, for $k\geq 2$, Ford \cite{ford} has shown that
there is a number $m$ such $\phi ^{-1}(m)$ contains precisely $k$ elements.\\

%% file: tot.bbl
\begin{thebibliography}{10}

\bibitem{chapman}
Chapman R., 
{\it Dirichlet's theorem: a real variable approach},
not published

\bibitem{ford} 
Ford K.,
{\it The number of solutions of $\phi (x)=m$},
Annals of mathematics, 150 (1999), 283-311.

\bibitem{gupta}
Gupta H.,
{\it Euler's totient function and its inverse},
Indian J. pure appl. Math., 12(1): 22-29, Jan 1981. 

 
\end{thebibliography}
